\documentclass[12pt,a4paper,reqno]{amsart}   

\usepackage{mathtools}  
\usepackage{ifthen} 
\newboolean{colorinverted}
\setboolean{colorinverted}{true} 
\setboolean{colorinverted}{false} 

\ifthenelse{\boolean{colorinverted}}{ 
\usepackage[usenames,dvipsnames,svgnames,table]{xcolor} 
\pagecolor{darkgray} 
\color{White} 
}

\usepackage{float}  
\usepackage{url}

\usepackage[myheadings]{fullpage}  
\usepackage[english]{babel} 
\usepackage{enumerate}   
\usepackage{datetime}

\usepackage{bm} 
\usepackage{mathrsfs}

\usepackage[full]{textcomp} 
\usepackage{newtxtext} 
\usepackage{cabin} 
\usepackage{zlmtt}
\usepackage[bigdelims,vvarbb]{newtxmath} 
\usepackage[cal=boondoxo]{mathalfa} 
 
\usepackage{microtype} %
  
\newcommand{\eps}{\varepsilon}

\newcommand{\Q}{\mathbb{Q}}

\renewcommand{\a}{\mathfrak{a}}

\renewcommand{\Re}{\textrm{Re}\,}

\newcommand{\odip}[2]{{o}_{#1}\left(#2\right)}
\newcommand{\odi}[1]{\odip{}{#1}}

\allowdisplaybreaks

\renewcommand{\qedsymbol}{$\square$}

\newtheoremstyle{sltheorems}
{10pt}
{6pt}
{\slshape}
{}
{\bfseries}
{.}
{.5em}
{\thmname{#1}\thmnumber{ #2}\thmnote{ (#3)}}

\theoremstyle{sltheorems} 
\newtheorem{Theorem}{Theorem}
\newtheorem{Corollary}{Corollary}
\newtheorem{Lemma}{Lemma}

\newtheoremstyle{remark}
{10pt}
{6pt}
{\rm} 
{}
{\bfseries}
{.}
{.5em}
{\thmname{#1}\thmnumber{ #2}\thmnote{ (#3)}}
 \theoremstyle{remark}

\usepackage{caption} 
\newcommand{\boundtheoretic}{404}
\newcommand{\boundcomp}{1000}
\newcommand{\boundtables}{200}
\allowdisplaybreaks

\begin{document} 
\let \oldsection \section
\def \newsection {\vspace{14pt plus 3pt}\oldsection}
\renewcommand{\section}{\newsection}

\title[Uniform effective estimates for  $\vert L(1,\chi)\vert$]
{Uniform effective estimates for  $\vert L(1,\chi)\vert$}  
 \author{Alessandro Languasco and Timothy S.~Trudgian}


\begin{abstract}   
Let $L(s,\chi)$ be the Dirichlet $L$-function associated to a
non-principal primitive Dirichlet character $\chi$ defined $\bmod\, q$, where $q\ge 3$.
We prove, under the assumption of the Generalised Riemann Hypothesis, 
 the validity of estimates given by Lamzouri, Li, and Soundararajan on 
$\vert L(1,\chi) \vert$. 
As a corollary, we  have that similar estimates hold for the class
number of the imaginary quadratic field  $\Q(\sqrt{-q})$, $q\ge 5$.
\end{abstract} 
\maketitle
\makeatletter
\def\subsubsection{\@startsection{subsubsection}{3}%
  \z@{.3\linespacing\@plus.5\linespacing}{-.5em}%
  {\normalfont\bfseries}} 
\makeatother
\mbox{}\vskip-1.5truecm
\section{Introduction}
Let $q\ge 3$ be an integer,  $\chi$ be a Dirichlet character $\bmod\, q$
and $L(s,\chi)$ be the associated Dirichlet $L$-function.
Assuming the Riemann Hypothesis for $L(s,\chi_d)$ holds,
where   $\chi_d$ is a quadratic Dirichlet character,
in 1928  Littlewood \cite{Littlewood1928} proved, for $d\ne m^2$, that
\begin{equation}
\label{Littlewood-bounds}
\Bigl(
\frac{12e^\gamma}{\pi^2}(1+\odi{1})\log \log  d
\Bigr)^{-1}
<
L(1,\chi_d)
<
2e^\gamma (1+\odi{1}) \log \log  d
\end{equation}
as $d$ tends to infinity, where $\gamma$ is the Euler--Mascheroni constant.
In 2015  
Lamzouri, Li, and Soundararajan \cite[Theorem 1.5]{LamzouriLS2015}
proved an
effective form of Littlewood's inequalities:
assuming the Generalised Riemann Hypothesis (GRH) holds, 
for every integer $q\ge 10^{10}$ and for every non-principal primitive character $\chi \bmod\, q$, 
they obtained that
\begin{equation}
\label{LLS-upper}
\vert L(1,\chi) \vert
\le 
2 e^\gamma \Bigl( \log \log q -\log 2 + \frac12 +\frac{1}{\log \log q} \Bigr)
\end{equation}
and
\begin{equation}
\label{LLS-lower}
\frac{1}{\vert L(1,\chi) \vert} 
\le 
\frac{12e^\gamma}{\pi^2}\Bigl( \log \log q - \log 2 + \frac12 +\frac{1}{\log \log q}  
+ \frac{14 \log \log q}{\log q}\Bigr).
\end{equation}
Such inequalities were numerically verified, for every prime $3\le q \le 10^7$, by Languasco \cite{Languasco2020}. We remark that the recent papers \cite{Bilu} and \cite{Ernvall} used (\ref{LLS-upper}) and (\ref{LLS-lower}) in the region $q\geq 10^{10}$. 
 
 The main goal of this paper is to enlarge the $q$-range of validity of \eqref{LLS-upper}-\eqref{LLS-lower} in the following
\begin{Theorem}
\label{Thm-main}
Assume the Generalised Riemann Hypothesis. Then, for $q\geq 3$ both \eqref{LLS-upper} and \eqref{LLS-lower}  hold 
for every non-principal primitive Dirichlet character $\bmod\, q$.
\end{Theorem}
The first step in the proof of Theorem \ref{Thm-main}
will be to extend the range of validity of the proof of Theorem 1.5
of \cite{LamzouriLS2015}
from $q\ge 10^{10}$ to $q\ge \boundtheoretic$.
Then using a computer program  
we will verify that  \eqref{LLS-upper}-\eqref{LLS-lower}  both hold 
for every non-principal  Dirichlet character $\bmod\, q$, where $q$ 
is composite, $3 \le q \le \boundcomp$. Combining this result
with the computations already performed for every prime
$3\le q \le 10^7$, by Languasco \cite{Languasco2020},
we obtained that  \eqref{LLS-upper}-\eqref{LLS-lower} hold for every integer $q\ge 3$.

Moreover, letting
\begin{equation}
\label{fg-def}
f(q):= \log \log q - \log 2 + 1/2 +1/\log \log q,
 \quad
g(q):= f(q) + 14 (\log \log q)/\log q,
\end{equation}
we numerically obtained that
 \begin{equation}
  \label{ineq-qsmall}
 \frac{\pi^2}{12 e^\gamma} \frac{3.25}{g(q)} < \vert L(1,\chi) \vert < 0.47 \cdot 2 e^\gamma   f(q)
  \end{equation}
 hold for $3 \le q \le \boundcomp$, 
 and, concerning Littlewood's estimates in \eqref{Littlewood-bounds}, that
 \begin{equation}
  \label{ineq-qsmall-ULILLI}
 0.425\cdot  \frac{\pi^2}{12 e^\gamma}  \log \log q < \vert L(1,\chi) \vert  <  0.71 \cdot 2 e^\gamma  \log \log q,
  \end{equation}
holds for $4\le q \le \boundcomp$. 
We note that Chowla \cite{Chowla} proved that there are 
infinitely many quadratic characters 
$\chi \bmod q$ such that 
$L(1, \chi) \geq (1 + o(1))e^{\gamma} \log\log q$. 
The coefficient of $\log\log q$ on the left side 
of \eqref{ineq-qsmall-ULILLI} is $0.1962\ldots$, 
which is a long way from the Chowla's 
$e^{\gamma} = 1.781\ldots$. 
We refer the reader to the excellent article by 
Granville and Soundararajan \cite{GS} for more on this topic.

Let now $q\ge 5$, $-q$ be a discriminant and $\chi_{-q}(n) = (-q \mid n)$
be the Kronecker symbol, which is a primitive character $\bmod\, q$.
Moreover, let $h(-q)$ denote the class number of the imaginary quadratic field $\Q(\sqrt{-q})$.
From Theorem \ref{Thm-main}  and the famous Dirichlet class number formula, \emph{i.e.},
\[
h(-q) = \frac{\sqrt{q}}{\pi} \, L(1,\chi_{-q}),
\] 
we obtain the following 
\begin{Corollary}
\label{cor-classnumber}
Assume the Generalised Riemann Hypothesis. Let
$-q$ be a discriminant  and let $h(-q)$  
denote the class number of $\Q(\sqrt{-q})$. For
every $q\ge 5$ we have that 
\[
h(-q)
\le 
\frac{2 e^\gamma}{\pi} \sqrt{q} \Bigl( \log \log q -\log 2 + \frac12 +\frac{1}{\log \log q} \Bigr)
\]
and
\[
h(-q)
\ge 
\frac{\pi}{12e^\gamma}\sqrt{q} \Bigl( \log \log q - \log 2 + \frac12 +\frac{1}{\log \log q}  
+ \frac{14 \log \log q}{\log q}\Bigr)^{-1}.
\] 
\end{Corollary}
Corollary \ref{cor-classnumber} extends the range of validity of Corollary 1.3 of
\cite{LamzouriLS2015} from $q\ge 10^{10}$ to $q\ge 5$. 

It is natural to ask whether similar results hold for bounds on $|\zeta(1+it)|$. This was noted in \cite[p.\ 2394]{LamzouriLS2015}, where bounds analogous to (\ref{LLS-upper}) and (\ref{LLS-lower}) were stated to hold for $t\geq 10^{10}$, under the assumption of the Riemann hypothesis. There are some additional computational difficulties in extending the work of this article to cover $\zeta(1+it)$. We leave this to future work.

\begin{figure}[ht] 
\includegraphics[scale=0.925,angle=0]{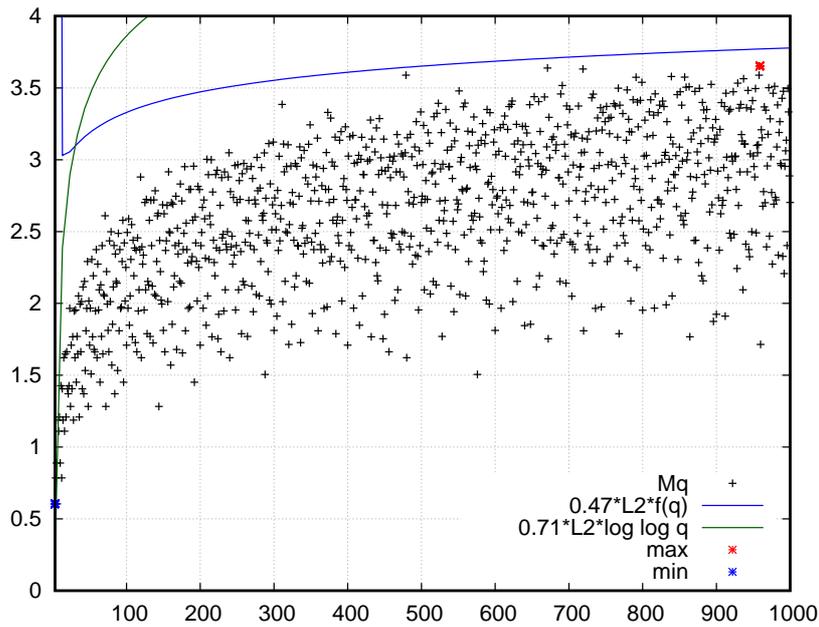}  
\caption{\small{The values of $M_q$,  $3\le q\le \boundcomp$.  
The minimal value for $M_q$ is $0.604599\dotsc$  attained at $q=3$ and the maximal one is  
$3.652103\dotsc$ attained at $q= 959$.    
The blue line represents $0.47\cdot L_2f(q)$, $L_2=  2 e^\gamma$,
where $f(q)$ is defined in \eqref{fg-def}. 
The green line represents $0.71\cdot L_2\log \log q$.  
 }}
\label{fig-LLS-1} 
 \end{figure}

\begin{figure}[ht] 
\includegraphics[scale=0.925,angle=0]{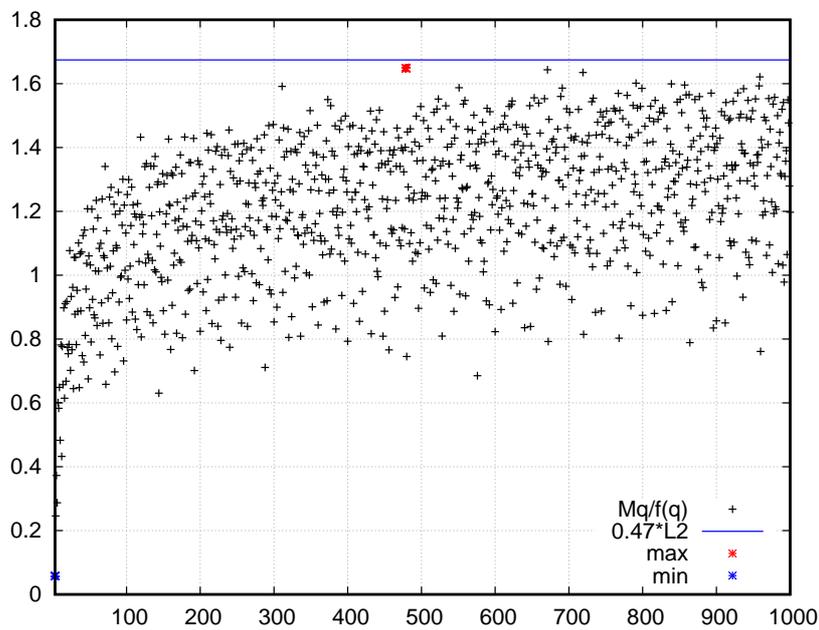}  
\caption{\small{The values of $M_q^\prime:=M_q/f(q)$, $3\le q\le \boundcomp$, 
where $f(q)$ is defined in \eqref{fg-def}. 
The minimal value for $M_q^\prime$ is $0.057396\dotsc$  attained at $q=3$ and the maximal one is  
$1.648945\dotsc$ attained at $q= 479$.  
The blue  line  represents $ 0.47\cdot L_2$,  $L_2=  2 e^\gamma$.
 }
 }
\label{fig-LLS-2} 
 \end{figure}

\begin{figure}[ht] 
\includegraphics[scale=1,angle=0]{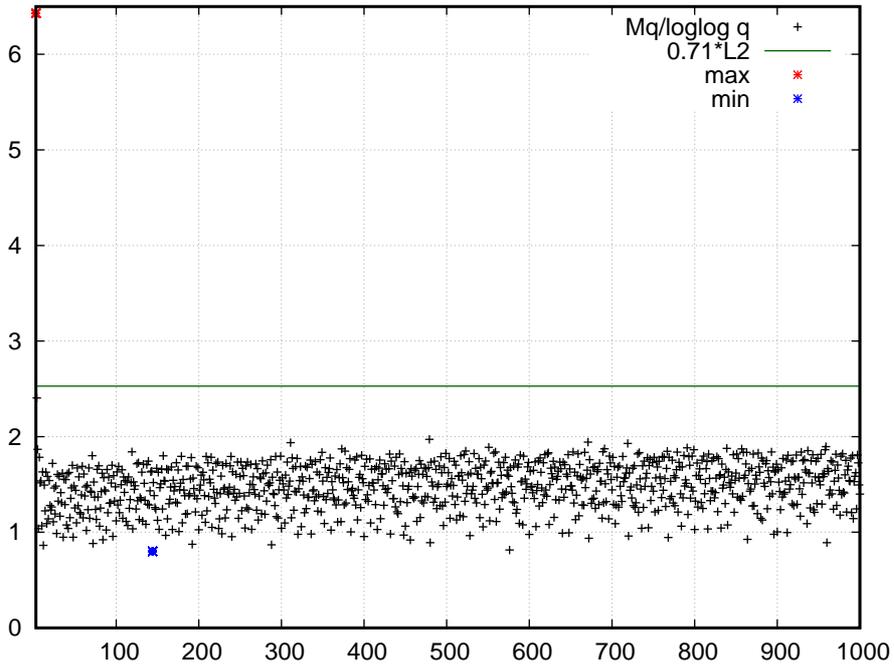}  
\caption{\small{The values of $M_q^{\prime\prime}:=M_q/\log \log q$, $3\le q\le \boundcomp$.
 The minimal value for $M_q^{\prime\prime}$ is $0.799902\dotsc$  attained at $q=144$ and the maximal one is   
$6.428641 \dotsc$ attained at $q= 3$.
The green  line  represents $0.71\cdot L_2$,  where  $L_2=  2 e^\gamma$. 
 }
 }
\label{fig-ULI} 
\end{figure}

\begin{figure}[ht] 
\includegraphics[scale=1,angle=0]{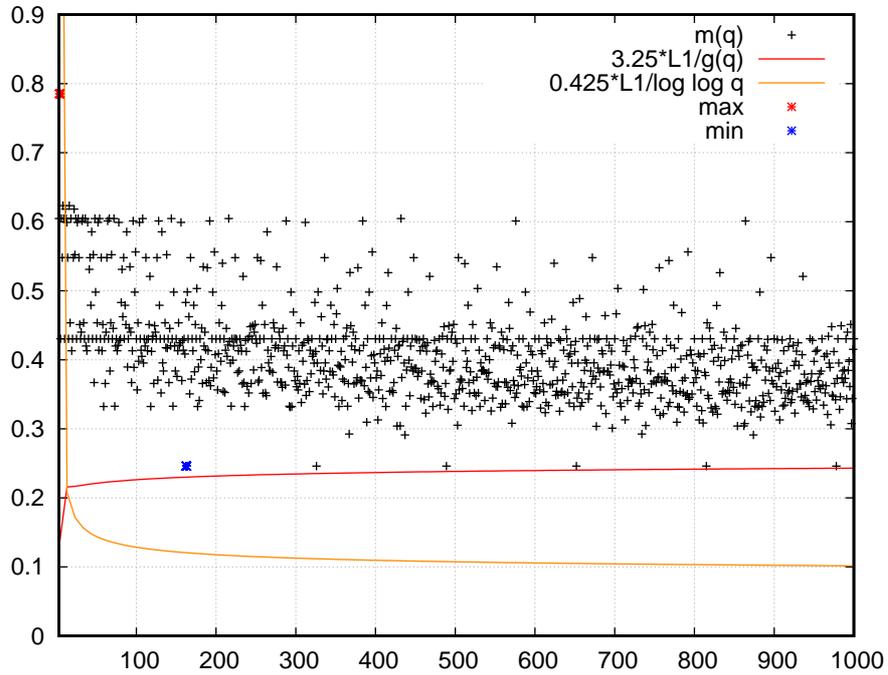}  
\caption{\small{ 
The values of $m_q$, $3\le q\le \boundcomp$.  
The minimal value for $m_q$ is $0.246068\dotsc$ attained at $q=163$
and the maximal one is   $0.785398\dotsc$ attained at $q=4$.
The red  line  represents  
$3.25 L_1 / g(q)$,  where $g(q)$ is defined in \eqref{fg-def}
and $L_1=   \frac{\pi^2}{12 e^\gamma}$. 
The orange line represents $0.425 L_1/\log\log q$. 
}
 }
\label{fig-LLS-3} 
\end{figure}

\begin{figure}[ht] 
\includegraphics[scale=1,angle=0]{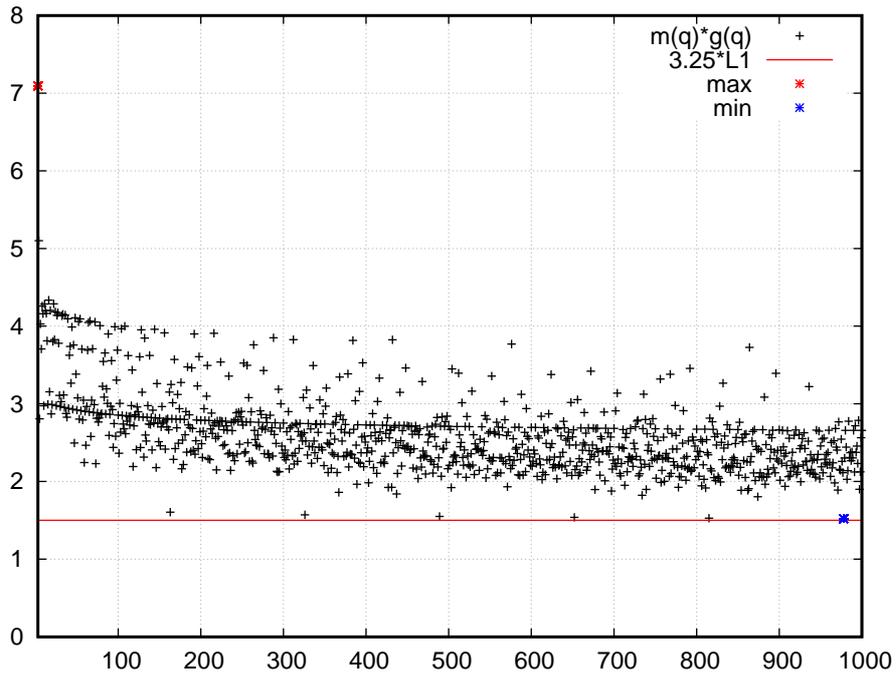}  
\caption{\small{The values of $m_q^\prime:=m_q g(q)$, $3\le q\le \boundcomp$, where $g(q)$ is defined in \eqref{fg-def}.  
The minimal value for $m_q^\prime$ is $1.520104\dotsc$ attained at $q=978$
and the maximal one is   $7.093329\dotsc$ attained at $q=3$.
The red  line  represents  
$3.25 L_1 $, where  $L_1=   \frac{\pi^2}{12 e^\gamma}$. 
 }
 }
\label{fig-LLS-4} 
 \end{figure}

\begin{figure}[ht] 
\includegraphics[scale=1,angle=0]{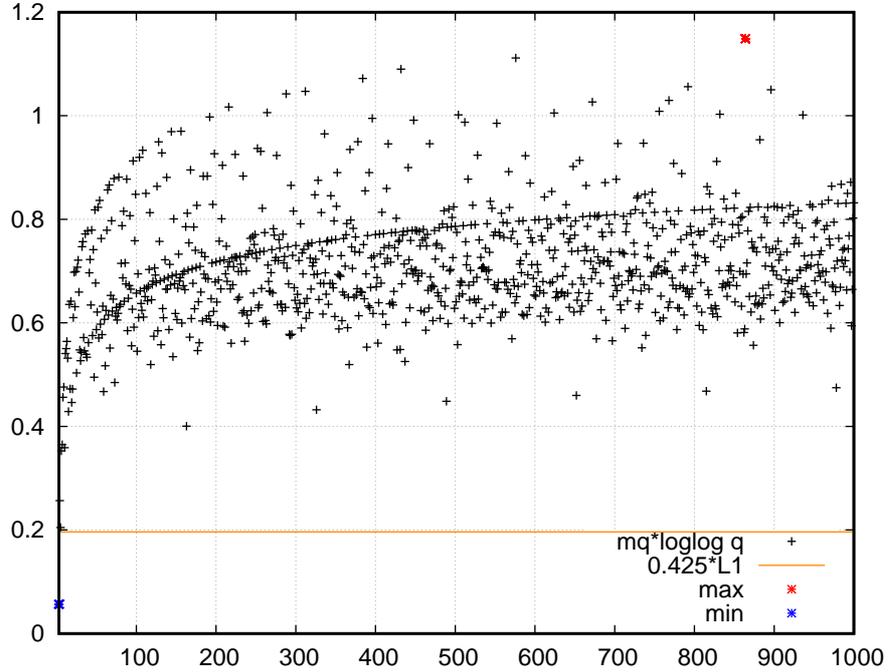}  
\caption{\small{The values of $m_q^{\prime\prime}:=m_q\log \log q$, 
$3\le q\le \boundcomp$ . 
The minimal value for $m_q^{\prime\prime}$ is $0.056861\dotsc$ attained at $q=3$ 
and the maximal one is   $1.148889\dotsc$ attained at $q=864$. 
The orange  line  represents  
$0.425 L_1$,  where   $L_1=   \frac{\pi^2}{12 e^\gamma}$.  
 }
 } 
\label{fig-LLI} 
\end{figure}

The paper is organised as follows:  in Section \ref{def-lemmas}
we state the needed lemmas from \cite{LamzouriLS2015}
and in Section \ref{main-thm-proof} we prove that 
\eqref{LLS-upper}-\eqref{LLS-lower} hold for every $q$ with $q\ge \boundtheoretic$.
Then in Section \ref{computL} we use a computer program
to verify the validity of \eqref{LLS-upper}-\eqref{LLS-lower}  for every $q$ with $3 \le q \le \boundcomp$.
In Section \ref{tables-figures}, after the references, we include Tables \ref{table1}-\ref{table2}
about the values  of $\vert L(1,\chi)\vert$, $\chi \bmod\, q$, $3 \le q \le \boundcomp$. These are augmented by
 Figures \ref{fig-LLS-1}-\ref{fig-LLI}, which contain several scatter plots. 
 
 \medskip 
\textbf{Acknowledgements}. TST is partially supported by ARC  DP160100932 and FT160100094.
 We wish to thank Youness Lamzouri for a discussion.
 \mbox{}\vskip-0.5truecm
\section{Definitions and Lemmas} 
\label{def-lemmas}
We require the following definitions and lemmas
from Lamzouri, Li, and Soundararajan \cite{LamzouriLS2015}.
Throughout the paper we assume the truth of the Generalised Riemann Hypothesis (GRH).

We let $\vartheta$ stand for a complex number of magnitude at most one. In each occurrence $\vartheta$
may stand for a different value, so that we may write $\vartheta-\vartheta = 2\vartheta$, $\vartheta \cdot \vartheta = \vartheta$
and so on. 
We recall that
\[
\xi(s) = s(s-1) \pi^{-s/2} \Gamma\Bigl(\frac{s}{2}\Bigr) \zeta(s),
\]
where $\zeta(s)$ is the Riemann zeta-function and $\Gamma(s)$
is Euler's function. The function $\xi(s)$ is an entire function of order $1$, satisfying
the functional equation $\xi(s)=\xi(1-s)$, and for which the Hadamard
factorisation formula gives
\[
\xi(s) = e^{Bs} \prod_\rho \Bigl( 1 - \frac{s}{\rho} \Bigr) e^{s/\rho},
\]
where $\rho$ runs over the non-trivial zeros of $\zeta(s)$, and $B$ is a real number 
given by
\begin{equation}
\label{B-def}
B = - \sum_\rho \Re  \frac{1}{\rho}  = \frac12 \log (4\pi) - 1 - \frac{\gamma} {2} = - 0.0230957089\dotsc.
\end{equation}
Let $L(s,\chi)$ be the Dirichlet $L$-function associated to a
non-principal primitive Dirichlet character $\chi$ defined $\bmod\, q$, where $q\ge 3$ 
is an integer. Let $\a=0$ if $\chi$ is even and $\a=1$ if $\chi$ is odd.
Let $\xi(s,\chi)$ be the completed $L$-function
\[
\xi(s,\chi)=\Bigl(\frac{q}{\pi}\Bigr)^{s/2} \Gamma\Bigl(\frac{s+\a}{2}\Bigr) L(s,\chi)
\]
which satisfies the functional equation
\[
\xi(s,\chi) = \eps_\chi  \xi(1-s,\overline{\chi}),
\]
where $\eps_\chi$ is a complex number of size $1$. The zeros of $\xi(s,\chi)$
are the non-trivial zeros of  $L(s,\chi)$, and letting $\rho_\chi =1/2+i \gamma_\chi$ denote
such a zero, we have Hadamard's factorisation formula:
\[
\xi(s,\chi)  = \exp(A(\chi) + sB(\chi)) \prod_{\rho_\chi} \Bigl( 1- \frac{s}{\rho_\chi} \Bigr) e^{s/\rho_\chi},
\]
where $A(\chi)$ and $B(\chi)$ are constants.  In particular
\[
\Re B(\chi)  = \Re B(\overline{\chi})  = \Re \frac{\xi^\prime}{\xi} (0, \chi) = - \sum_{\rho_\chi}\Re \frac{1}{\rho_\chi}.
\]
Recall now the digamma function $\psi(z) =\Gamma^\prime/\Gamma(z)$;
we will need the following special values
\[
\psi(1) = -\gamma,\quad
\psi\Bigl(\frac12\Bigr) = -2\log2 -\gamma.
\]

Using such definitions, we now state the lemmas we need from \cite{LamzouriLS2015}.
\begin{Lemma}[Lemma 2.3 of \cite{LamzouriLS2015}]
\label{LLS-lemma-2.3}
Assume GRH. Let $q\ge 3$ and let $\chi \bmod\, q$ be a primitive Dirichlet character.
For any $x>1$, we have, for some $\vert \vartheta \vert \le 1$, that
\begin{align*}
- \frac{\xi^\prime}{\xi} (0, \overline{\chi})
- \frac{1}{x}\, \frac{\xi^\prime}{\xi} (0, \chi)
+ \frac{2\vartheta}{\sqrt{x}} \vert \Re B(\chi) \vert
=
\frac12 \Bigl(1-\frac1x \Bigr) \log \frac{q}{\pi}
-
\sum_{n\le x} \frac{\Lambda(n)\chi(n)}{n}  \Bigl(1-\frac{n}{x} \Bigr)
+ E_{\a}(x),
 \end{align*}
 where
 \begin{equation}
 \label{E0-def}
E_0(x) := -\log 2 - \frac{\gamma}{2} \Bigl(1-\frac1x \Bigr) + \frac{\log x +1}{x} 
- \sum_{k=1}^{+\infty} \frac{x^{-2k-1}}{2k(2k+1)}
\end{equation}
and
 \begin{equation}
 \label{E1-def}
E_1(x) := 
- \sum_{k=0}^{+\infty} \frac{x^{-2k-2}}{(2k+1)(2k+2)}
 - \frac{\gamma}{2} \Bigl(1-\frac1x \Bigr) + \frac{\log 2}{x} .
\end{equation}
In particular,  $\vert \Re B(\chi) \vert$ equals
\[
\Bigl(1 + \frac{2\vartheta}{\sqrt{x}} + \frac1x \Bigr)^{-1}
\Bigl(\frac12 \Bigl(1-\frac1x \Bigr)  \log \frac{q}{\pi}  
-     \Re \sum_{n\le x} \frac{\Lambda(n)\chi(n)}{n}  \Bigl(1-\frac{n}{x} \Bigr)
+ E_{\a}(x) \Bigr).
\]
\end{Lemma}

\begin{Lemma}[Lemma 2.4 of \cite{LamzouriLS2015}]
\label{LLS-lemma-2.4}
Assume RH. For $x>1$ we have, for some $\vert \vartheta \vert \le 1$, that
\[
\Sigma_2
:=
 \sum_{n\le x} \frac{\Lambda(n)}{n}  \Bigl(1-\frac{n}{x} \Bigr)
 = 
 \log x -(1+\gamma)+ \frac{\log (2\pi)}{x} 
 -
 \sum_{n=1}^{+\infty} \frac{x^{-2n-1}}{2n(2n+1)}
 + 2\vartheta \frac{B}{\sqrt{x}}.
\]
\end{Lemma}

\begin{Lemma}[Lemma 2.5 of \cite{LamzouriLS2015}]
\label{LLS-lemma-2.5}
Let $q\ge 3$ and let $\chi \bmod\, q$ be a primitive Dirichlet character.
Suppose that GRH holds for $L(s,\chi)$. For any $x\ge 2$, there exists
a real number $\vert \vartheta \vert \le 1$ such that
\begin{align*}
\log \vert L(1,\chi) \vert 
&= 
\Re \sum_{n\le x} \frac{\chi(n) \Lambda(n)}{n \log n} \frac{\log (x/n)}{\log x}
+
\frac{1}{2\log x}\Bigl(  \log \frac{q}{\pi} +   \psi\Bigl(\frac{1+\a}{2}\Bigr)\Bigr)
\\ &\hskip1cm
-
\Bigl(\frac{1}{\log x}  + \frac{2\vartheta}{\sqrt{x}(\log x)^2}\Bigr) \vert \Re B(\chi)\vert 
+ 
\frac{2 \vartheta}{x(\log x)^2}.
 \end{align*}
\end{Lemma}

\begin{Lemma}[Lemma 2.6 of \cite{LamzouriLS2015}]
\label{LLS-lemma-2.6}
Assume RH. For all $x\ge e$, there exists
a real number $\vert \vartheta \vert \le 1$ such that 
\[
\Sigma_1
:=
 \sum_{n\le x} \frac{\Lambda(n)}{n \log n}  \frac{\log (x/n)}{\log x} 
 = 
\log \log x  +\gamma -1 + \frac{\gamma}{\log x}  
+ \frac{2B\vartheta}{\sqrt{x}(\log x)^2}
+ \frac{\vartheta}{3 x^3(\log x)^2}
 .
\]
\end{Lemma}

\begin{Lemma}[Lemma 5.1 of \cite{LamzouriLS2015}]
\label{LLS-lemma-5.1}
Let $q\ge 3$ and let $\chi \bmod\, q$ be a primitive Dirichlet character.
For $x\ge 3$  we have
\begin{equation}\label{verdi}
\Re
\sum_{n\le x }\Lambda(n) \chi(n) \Bigl(  \frac{1}{n \log n}  - \frac{1}{x\log x} \Bigr)
\ge
\sum_{p^k\le x }\Lambda(p^k)  (-1)^k \Bigl(  \frac{1}{p^k \log p^k}  - \frac{1}{x\log x} \Bigr).
\end{equation}
\end{Lemma}
We remark that Lemma \ref{LLS-lemma-5.1} was stated just for $x\ge 100$ in \cite{LamzouriLS2015},
but it is easy to verify that its proof holds for $x\ge 3$.
\mbox{}\vskip-0.5truecm
\section{Proof of Theorem \ref{Thm-main}}
\label{main-thm-proof} 

From now on we assume $3\le q < 10^{10}$ since 
the result was already proved for $q\ge 10^{10}$ in \cite{LamzouriLS2015}.

 \subsection{Upper bound} 
 \label{upperbound}
Let $x \ge 2$. Using Lemma \ref{LLS-lemma-2.5} and
recalling that  $\vert \vartheta \vert \le 1$  we have
\begin{align}\label{agrippina}
\log \vert L(1,\chi) \vert 
&\le 
\Re \sum_{n\le x} \frac{\chi(n) \Lambda(n)}{n \log n} \frac{\log (x/n)}{\log x}
+
\frac{1}{2 \log x}\Bigl(  \log \frac{q}{\pi} +   \psi\Bigl(\frac{1+\a}{2}\Bigr)\Bigr)\notag
\\ &\hskip1cm
-
\Bigl(\frac{1}{\log x}  - \frac{2}{\sqrt{x}(\log x)^2}\Bigr) \vert \Re B(\chi)\vert 
+ 
\frac{2}{x(\log x)^2}.
 \end{align}
 For $x\ge 2$, we can use Lemma \ref{LLS-lemma-2.3}
 so that
 \begin{align*}
 \vert \Re B(\chi) \vert
 \ge
\Bigl(1 + \frac{1}{\sqrt{x}} \Bigr)^{-2}
\Bigl(\frac12 \Bigl(1-\frac1x \Bigr)  \log \frac{q}{\pi}  
-     \Re \sum_{n\le x} \frac{\Lambda(n)\chi(n)}{n}  \Bigl(1-\frac{n}{x} \Bigr)
+ E_{\a}(x) \Bigr).
 \end{align*} 
 
We wish to show, for $\a\in\{0,1\}$ and  $x\geq 9$, that
\begin{equation}
\label{whisky}
C_{\a}^{+}(x)
:= 
-E_{\a}(x) 
\Bigl(1 + \frac{1}{\sqrt{x}} \Bigr)^{-2} 
\Bigl( \frac{1}{\log x} - \frac{2}{\sqrt{x}(\log x)^2}\Bigr) 
+ \frac{1}{2 \log x} \psi \Bigl( \frac{1 + \a}{2}\Bigr) 
+ \frac{2}{x (\log x)^2}
\le 0.
\end{equation}
This is equivalent to proving that $D^{+}_{\a}(x) := C^{+}_{\a}(x) \log x \le 0$
for $\a\in\{0,1\}$ and  $x\geq 9$.

We first study the case $\a=0$ of \eqref{whisky}. We begin by proving that $E_{0}(x) <0$ for $x\geq 2$.
To do this, let us consider upper and lower bounds on the sum over $k$ in $E_{0}(x)$ 
in \eqref{E0-def}. We have, for any $K\geq 1$,  that
\begin{align}
\notag
\sum_{k=1}^{\infty} \frac{x^{-2k-1}}{2k(2k+1)} &= \sum_{k=1}^{K} \frac{x^{-2k-1}}{2k(2k+1)} + \sum_{k\geq K+1} \frac{x^{-2k-1}}{2k(2k+1)}\\
\notag
&\leq \sum_{k=1}^{K} \frac{x^{-2k-1}}{2k(2k+1)} + \frac{1}{x^{3+2K}} \sum_{k\geq K+1} \frac{1}{2k(2k+1)}\\
\notag
&\leq \sum_{k=1}^{K} \frac{x^{-2k-1}}{2k(2k+1)} + \frac{1}{4x^{3+2K}} \sum_{k\geq K+1} \frac{1}{k^{2}}\\
&
\label{bourbon} 
\leq \sum_{k=1}^{K} \frac{x^{-2k-1}}{2k(2k+1)} + \frac{1}{4K x^{3+ 2K}}.
\end{align}
Clearly, for any $K\geq 1$, we also have
\[
 \sum_{k=1}^{\infty} \frac{x^{-2k-1}}{2k(2k+1)} \geq  \sum_{k=1}^{K} \frac{x^{-2k-1}}{2k(2k+1)}.
 \]
Hence
\begin{equation}\label{sake}
- E_{0}(x) \geq 
\log 2 + \frac{\gamma}{2} \Bigl(1-\frac1x \Bigr)
- \frac{\log x +1}{x} 
+ \sum_{k=1}^{K} 
\frac{ x^{-2k-1}}{2k(2k+1)}.
\end{equation}
Taking $K=2$ we see that right-hand side in \eqref{sake} is positive for all $x\geq 2$.
Hence, remarking that $ (1 + 1/\sqrt{x} )^{-2}  \le 1$ and that 
$1- 2/(\sqrt{x}\log x) >0$ for every $x \ge 4$, we have that
\begin{equation}
\label{gin}
- E_0(x) 
\le 
F_0^{+}(x) 
:=
\log 2 + \frac{\gamma}{2} \Bigl(1-\frac1x \Bigr)
- \frac{\log x +1}{x} 
+ \sum_{k=1}^{K} 
\frac{ x^{-2k-1}}{2k(2k+1)}
+  \frac{1}{4K x^{3+ 2K}}
\end{equation}
 and that
\begin{align}
\notag
D_0^{+}(x) 
&\le  
 \frac{2}{x\log x} - \frac{\gamma}{2x}  
- \frac{\log x +1}{x} 
+ \sum_{k=1}^{K}  \frac{ x^{-2k-1}}{2k(2k+1)}
+  \frac{1}{4K x^{3+ 2K}}
- \frac{2}{\sqrt{x} \log x}F_0^{+}(x)
\\
\label{vodka}
&
= G_0^{+}(x) - \frac{2}{\sqrt{x} \log x}F_0^{+}(x),
\end{align}
say. Choosing  $K=2$ in \eqref{gin} and \eqref{vodka},
 we immediately have that $G_0^{+}(x)<0$ 
and $F_0^{+}(x) > 0$ for every $x\ge 3$.
This proves that $D_0^{+}(x)  <0$  for every $x\ge 4$. Hence \eqref{whisky} holds for $\a=0$ and $x\ge 4$.

We now study the case $\a=1$ of \eqref{whisky}. 
This is done in the same way, noting that $\psi(1) = -\gamma$  in  \eqref{whisky}
accounts for the main term in $E_{1}(x)$ in \eqref{E1-def}. 
We estimate the sum over $k$ in the definition of $E_{1}(x)$ in 
\eqref{E1-def} in the same way as in \eqref{bourbon}. This gives
\begin{equation}
\label{brandy}
\sum_{k=0}^{\infty} \frac{x^{-2k-2}}{(2k+1)(2k+2)} 
\leq 
\sum_{k=0}^{K} \frac{x^{-2k-2}}{(2k+1)(2k+2)} 
+ 
\frac{1}{4K x^{4 + 2K}}.
\end{equation}
Using this we see that $E_{1}(x)\leq 0$ for all $x\geq 4$. 
We therefore use $(1 + 1/\sqrt{x})^{-2} \leq 1$ and have
\begin{equation}
\label{cognac}
-E_{1}(x) 
\leq 
F_{1}^{+}(x)
:= 
\frac{\gamma}{2} \Bigl( 1 - \frac{1}{x}\Bigr) 
- \frac{\log 2}{x} 
+ \sum_{k=0}^{K} \frac{x^{-2k-2}}{(2k+1)(2k+2)} 
+ \frac{1}{4K x^{4 + 2K}}.
\end{equation}
We proceed to obtain a bound similar to \eqref{vodka}, namely
\begin{align}
\notag
D_1^{+}(x) &
\leq 
- \frac{\gamma}{2x} 
- \frac{\log 2}{x} 
+ \frac{2}{x \log x} 
+ \sum_{k=0}^{K} \frac{x^{-2k-2}}{(2k+1)(2k+2)} 
+ \frac{1}{4K x^{2K+4}} 
- \frac{2}{\sqrt{x} \log x} F_{1}(x)
\\
\label{amaretto} 
&
= 
G_{1}^{+}(x) - \frac{2}{\sqrt{x} \log x} F_{1}^{+}(x),
\end{align} 
say. Taking $K=1$ we see that $G_{1}^{+}(x)<0$ for all $x\geq 9$. 
Since we have $F_{1}^{+}(x) >0$ for $x\geq 4$, we conclude that $D_1^{+}(x)\leq 0$ for $x\geq 9$.
Hence \eqref{whisky} holds for $\a=1$ and $x\ge 9$.

Combining the previous results we can write that \eqref{whisky} holds for $\a\in\{0,1\}$ and $x\ge 9$.

We now combine all our results, keeping in place the minor terms that, for large $x$ were simplified in \cite{LamzouriLS2015}. 
We find that
\begin{align}
\notag
    \log \vert L(1, \chi) \vert &\leq \Sigma_{1} 
    + \frac{\log (q/\pi)}{2\log x} \Bigl\{ 1 - \Bigl( 1 - \frac{1}{x}\Bigr)
    \Bigl( 1 - \frac{2}{\sqrt{x} \log x}\Bigr)
    \Bigl(1 + \frac{1}{\sqrt{x}}\Bigr)^{-2} \Bigr\} 
    \\&
    \label{lagavulin}
    + \frac{1}{\log x} \Bigl( 1 - \frac{2}{\sqrt{x} \log x}\Bigr)
    \Bigl( 1 + \frac{1}{\sqrt{x}}\Bigr)^{-2} \Sigma_{2},
\end{align}
where we use Lemma \ref{LLS-lemma-2.6} to bound $\Sigma_{1}$ and 
Lemma \ref{LLS-lemma-2.4} to bound $\Sigma_{2}$, taking just one 
term in the sum over $n$. Putting all of this together, we set 
$x = 1/4(\log q)^{2}$, as in \cite{LamzouriLS2015}, and exponentiate. 
We aim at showing that the resultant expression does not exceed 
the bound in \eqref{LLS-upper}, subject to the restriction that 
$x\geq 9$. We have that \eqref{lagavulin} implies \eqref{LLS-upper}
whenever $q\geq 15$, however, the bound $x\geq 9$ means that 
we only win for $q\geq \boundtheoretic$. 

 \subsection{Lower bound} 
Let $x \ge 2$. We begin with a lower bound on $\log |L(1, \chi)|$ similar to that obtained in (\ref{agrippina}). Again for $x\ge 2$, we can use Lemma \ref{LLS-lemma-2.3}
 so that
 \begin{align*}
 \vert \Re B(\chi) \vert
 \le
\Bigl(1 - \frac{1}{\sqrt{x}} \Bigr)^{-2}
\Bigl(\frac12 \Bigl(1-\frac1x \Bigr)  \log \frac{q}{\pi}  
-     \Re \sum_{n\le x} \frac{\Lambda(n)\chi(n)}{n}  \Bigl(1-\frac{n}{x} \Bigr)
+ E_{\a}(x) \Bigr).
 \end{align*} 
 
 We wish to show, for $\a\in\{0,1\}$ and  $x\geq 5$, that
\begin{equation}
\label{ale}
C_{\a}^{-}(x)
:= 
-E_{\a}(x) 
\Bigl(1 - \frac{1}{\sqrt{x}} \Bigr)^{-2} 
\Bigl( \frac{1}{\log x} + \frac{2}{\sqrt{x}(\log x)^2}\Bigr) 
+ \frac{1}{2 \log x} \psi \Bigl( \frac{1 + \a}{2}\Bigr) 
- \frac{2}{x (\log x)^2}
\ge 0.
\end{equation}
This is equivalent to proving that $D^{-}_{\a}(x) := C^{-}_{\a}(x) \log x \ge 0$
for $\a\in\{0,1\}$ and  $x\geq 9$. We recall that 
in \S \ref{upperbound} we proved that $E_0(x)<0$ for $x\ge 2$ and $E_1(x)<0$ for $x\ge 4$.

We first study the case $\a=0$ of \eqref{ale}. We use $(1 - 1/\sqrt{x})^{-2} \geq 1 + 2/\sqrt{x}$ 
and take $K=1$ to see that $C_{0}^{-}(x)\geq 0$ for $x\geq 6$.

We now study the case $\a=1$ of \eqref{ale}. 
Again, we need to use the inequality  $(1 - 1/\sqrt{x})^{-2} \geq 1 + 2/\sqrt{x}$ 
rather than just the trivial bound of $\geq 1$. Taking $K=1$ shows that $C_{1}^{-}(x)\geq 0$ for $x\geq 9$.

Now we proceed to prove the lower bound.
Keeping all terms together we invoke
Lemmas \ref{LLS-lemma-2.4} and \ref{LLS-lemma-2.6}, 
and only require $x\geq 9$ since we just proved that then $C_{\a}^{-}(x) \geq 0$ for $\a\in\{0,1\}$. 
Note that Lemma \ref{LLS-lemma-5.1} holds for any $x\ge 3$. 
Note also that we only require lower bounds on $\sum_{n\geq 1} x^{-2n-1}/(2n(2n+1))$, 
whence we just take the term $n=1$. 
Arguing as on page 2408 
of \cite{LamzouriLS2015} and
putting everything together we arrive at the following, for $x\geq 9$,
\begin{align}
\notag
\log \vert L(1, \chi) \vert \geq
& -\Bigl(\log\log x + \gamma - 1 + \frac{\gamma}{\log x} + \frac{2 \vert  B \vert}{\sqrt{x} \log^{2} x}\Bigr)
\\& \notag
- \frac{1}{\log x} \Bigl( \log x - (1 + \gamma) + \frac{\log 2\pi}{x} - \frac{1}{6x^{3}} + \frac{2 \vert B \vert}{\sqrt{x}}\Bigr) 
+ \log \zeta(2) - \frac{3}{2 \sqrt{x}}
\\&\notag
-\Bigl\{ \Bigl( 1 - \frac{1}{\sqrt{x}}\Bigr)^{-2} 
\Bigl( 1 + \frac{2}{\sqrt{x} \log x}\Bigr) - 1\Bigr\} 
\Bigl( -1 + \frac{1 + \gamma}{\log x} - \frac{a(x)}{\log x}\Bigr)
\\ & 
\label{champagne} 
+ \frac{ \log (q/\pi)}{2 \log x} 
\Bigl\{ 1 - \Bigl( 1 - \frac{1}{x}\Bigr) \Bigl( 1 - \frac{1}{\sqrt{x}}\Bigr)^{-2} 
\Bigl( 1 + \frac{2}{\sqrt{x} \log x}\Bigr)\Bigr\},
\end{align}
where $a(x)$ comes from Lemma \ref{LLS-lemma-2.4}, and is defined below. Note that the $\log\zeta(2)$ term above comes from the contribution of the $k=2$ terms in the sum on the right of (\ref{verdi}). Recall that $B = -0.0230\ldots$ by \eqref{B-def}.

We have, from Lemma \ref{LLS-lemma-2.4},  
that
\begin{align}
\notag
\sum_{n\leq x} \frac{\Lambda(n)}{n} 
\Bigl( 1 - \frac{n}{x}\Bigr) &
= \log x - (1 + \gamma) + \frac{\log 2\pi}{x} 
- \sum_{n=1}^{\infty} \frac{x^{-2n-1}}{2n(2n+1)} + \frac{2 \vartheta B}{\sqrt{x}}
\\
\notag
& \leq \log x - (1+ \gamma) + \frac{\log 2\pi}{x} - \frac{1}{6x^{3}} - \frac{2 B}{\sqrt{x}}
\\
\label{IPA} 
& := \log x - (1+ \gamma) + a(x).
\end{align}
Note then, that $a(x)$ is positive and tends to zero as 
$x \to \infty$.
We now use the definition of $a(x)$ from \eqref{IPA} in 
\eqref{champagne}, 
and exponentiate, taking $x = 1/4  (\log q)^{2}$.  
We find that the bound obtained by \eqref{champagne} 
implies \eqref{LLS-lower} whenever $q\geq 120$. 
However, as before, we must satisfy the condition 
that $x\geq 9$. With our choice of $x = 1/4 (\log q)^{2}$ 
we may only deduce the result in the 
theorem when $q\geq \boundtheoretic$. 
%
\subsection{Computation of $\vert L(1,\chi)\vert$ for $3\le q \le \boundcomp$} 
 \label{computL}
Letting
\begin{equation}
\label{Max-def}
M_q:=\max_{\chi \neq \chi_0} 
\vert L(1,\chi) \vert 
\quad
\textrm{and}
\quad
m_q:=\min_{\chi \neq \chi_0} 
\vert L(1,\chi) \vert, 
\end{equation} 
we will obtain $M_q, m_q$  using PARI/GP, v.~2.13.0, 
since it has  the ability to  generate
 Dirichlet $L$-functions (and  many other $L$-functions).
 This can be done
with few instructions of the gp scripting language.
Such a computation has a linear cost in the number of calls
of the {\tt lfun} function of PARI/GP. 
We were able to get the values of $M_q, m_q$ for every integer $q$, $3\le q \le \boundcomp$, 
with a precision of  $30$ decimal digits within $19$ minutes and $57$ seconds
of computation time (or $13$ minutes and $4$ seconds for just the composite integers we needed).
In Tables \ref{table1} and \ref{table2} we present such data for $3\le q \le \boundtables$; to get them we needed
 less than  24 seconds of computation time for each table. 
The machine we used was a Dell OptiPlex-3050, equipped with an Intel i5-7500 processor, 3.40GHz, 
16 GB of RAM and running Ubuntu 18.04.4.
The analysis on these data to verify the inequalities in \eqref{ineq-qsmall}-\eqref{ineq-qsmall-ULILLI}
were performed using a python-pandas program. All the programs used here
can be downloaded from the following address:
\url{https://www.math.unipd.it/~languasc/LLS_ineq.html}.
    
 \mbox{}\vskip-.5truecm
\renewcommand{\bibliofont}{\normalsize} 

\vskip 0.5cm
\noindent
Alessandro Languasco,
Universit\`a di Padova,
 Dipartimento di Matematica,
 ``Tullio Levi-Civita'',
Via Trieste 63,
35121 Padova, Italy.   
{\it e-mail}: alessandro.languasco@unipd.it 

\medskip
\noindent
Timothy S.~Trudgian,   
School of Science, The University of New South Wales Canberra at ADFA,
ACT, 2610 Australia.
{\it e-mail}:  t.trudgian@adfa.edu.au\\

\newpage 

\section{Tables and figures}
\label{tables-figures}
\begin{table}[H]
\scalebox{0.606}{
\begin{tabular}{|c|c|}
\hline
$q$  &  $M_q$\\ \hline 
$3$ & $0.60459978807807261686469275254$ \\
$4$ & $0.78539816339744830961566084582$ \\
$5$ & $0.88857658763167324940317619801$ \\
$6$ & $0.60459978807807261686469275254$ \\
$7$ & $1.18741041172372594878462529795$ \\
$8$ & $1.11072073453959156175397024752$ \\
$9$ & $1.20919957615614523372938550509$ \\
$10$ & $0.88857658763167324940317619801$ \\
$11$ & $1.42640418224108352050983157388$ \\
$12$ & $0.78539816339744830961566084582$ \\
$13$ & $1.40613477980703732992641904009$ \\
$14$ & $1.18741041172372594878462529795$ \\
$15$ & $1.62231147038944475878118430812$ \\
$16$ & $1.11072073453959156175397024752$ \\
$17$ & $1.64849370699838393605712538480$ \\
$18$ & $1.20919957615614523372938550509$ \\
$19$ & $1.66331503401599345646761861198$ \\
$20$ & $1.40496294620814527863127492864$ \\
$21$ & $1.37110344169451507464463658377$ \\
$22$ & $1.42640418224108352050983157388$ \\
$23$ & $1.96520205410785916590276700512$ \\
$24$ & $1.28254983016186409554403635967$ \\
$25$ & $1.70951180694329794837835134968$ \\
$26$ & $1.40613477980703732992641904009$ \\
$27$ & $1.76773152072033877884666136808$ \\
$28$ & $1.18741041172372594878462529795$ \\
$29$ & $1.94760835256298812067787472221$ \\
$30$ & $1.62231147038944475878118430812$ \\
$31$ & $1.95460789577555396208826644691$ \\
$32$ & $1.45122657606971550702884393708$ \\
$33$ & $1.76974436700717745842473969267$ \\
$34$ & $1.64849370699838393605712538480$ \\
$35$ & $2.05172721872421812035646128332$ \\
$36$ & $1.20919957615614523372938550509$ \\
$37$ & $1.99550481523309639952094476513$ \\
$38$ & $1.66331503401599345646761861198$ \\
$39$ & $2.01222972650783293991026647857$ \\
$40$ & $1.40496294620814527863127492864$ \\
$41$ & $2.11431182971789691740711719452$ \\
$42$ & $1.37110344169451507464463658377$ \\
$43$ & $2.15300367350566872872181166589$ \\
$44$ & $1.53264358003262663869763307467$ \\
$45$ & $1.80945332629550203732956607281$ \\
$46$ & $1.96520205410785916590276700512$ \\
$47$ & $2.29124192852861593669991478644$ \\
$48$ & $1.28254983016186409554403635967$ \\
$49$ & $2.01688417594425696685721289539$ \\
$50$ & $1.70951180694329794837835134968$ \\
$51$ & $1.96734163237258188833523337933$ \\
$52$ & $1.50917229488087470493269985892$ \\
$53$ & $2.30607194293623581454960238070$ \\
$54$ & $1.76773152072033877884666136808$ \\
$55$ & $2.19228536898779746002203139574$ \\
$56$ & $1.67925190836271397945223134714$ \\
$57$ & $1.66331503401599345646761861198$ \\
$58$ & $1.94760835256298812067787472221$ \\
$59$ & $2.37483432765382947109503698217$ \\
$60$ & $1.62231147038944475878118430812$ \\
$61$ & $2.27383316907813941451707465286$ \\
$62$ & $1.95460789577555396208826644691$ \\
$63$ & $2.09439510239319549230842892219$ \\
$64$ & $1.45122657606971550702884393708$ \\
$65$ & $2.12917795839603170191079901461$ \\
$66$ & $1.76974436700717745842473969267$ \\
$67$ & $2.40177874951129444901523775053$ \\
$68$ & $1.64849370699838393605712538480$ \\
\hline
\end{tabular}
}
\scalebox{0.606}{
\begin{tabular}{|c|c|}
\hline
$q$  &  $M_q$\\ \hline
 $69$ & $2.06213746115420018933849107128$ \\
$70$ & $2.05172721872421812035646128332$ \\
$71$ & $2.60986917715784586434512887899$ \\
$72$ & $1.28254983016186409554403635967$ \\
$73$ & $2.22099352575696724500691589494$ \\
$74$ & $1.99550481523309639952094476513$ \\
$75$ & $1.81177664489821466934352802615$ \\
$76$ & $1.66331503401599345646761861198$ \\
$77$ & $2.43523343388619715784941500753$ \\
$78$ & $2.01222972650783293991026647857$ \\
$79$ & $2.49865862902662621045751915333$ \\
$80$ & $1.57079632679489661923132169164$ \\
$81$ & $1.85460990010218743941830369297$ \\
$82$ & $2.11431182971789691740711719452$ \\
$83$ & $2.40279523907172221214735958325$ \\
$84$ & $1.37110344169451507464463658377$ \\
$85$ & $2.34387207448910709235492592796$ \\
$86$ & $2.15300367350566872872181166589$ \\
$87$ & $2.14059432707202598850826022652$ \\
$88$ & $1.53264358003262663869763307467$ \\
$89$ & $2.48752834330666367191078756009$ \\
$90$ & $1.80945332629550203732956607281$ \\
$91$ & $2.37482082344745189756925059590$ \\
$92$ & $1.96520205410785916590276700512$ \\
$93$ & $1.97538604461783624686153124312$ \\
$94$ & $2.29124192852861593669991478644$ \\
$95$ & $2.57856484292230594522273621834$ \\
$96$ & $1.45122657606971550702884393708$ \\
$97$ & $2.42085614235917869433795064647$ \\
$98$ & $2.01688417594425696685721289539$ \\
$99$ & $2.25006635803238006730577318676$ \\
$100$ & $1.70951180694329794837835134968$ \\
$101$ & $2.49348309598992905601857403054$ \\
$102$ & $1.96734163237258188833523337933$ \\
$103$ & $2.58872219777793564220546520282$ \\
$104$ & $1.84835102820162619473837967628$ \\
$105$ & $2.05172721872421812035646128332$ \\
$106$ & $2.30607194293623581454960238070$ \\
$107$ & $2.54845309686851323582504835532$ \\
$108$ & $1.76773152072033877884666136808$ \\
$109$ & $2.43977808110771365276355432077$ \\
$110$ & $2.19228536898779746002203139574$ \\
$111$ & $2.38549422924165041465092992869$ \\
$112$ & $1.72727422084656605141455939338$ \\
$113$ & $2.37003547987807971651823457293$ \\
$114$ & $1.66331503401599345646761861198$ \\
$115$ & $2.45609167047081871331670973715$ \\
$116$ & $1.94760835256298812067787472221$ \\
$117$ & $2.09439510239319549230842892219$ \\
$118$ & $2.37483432765382947109503698217$ \\
$119$ & $2.87989326382063351362429041213$ \\
$120$ & $1.62231147038944475878118430812$ \\
$121$ & $2.39329849251482613435551192673$ \\
$122$ & $2.27383316907813941451707465286$ \\
$123$ & $2.33588669690223816457514930983$ \\
$124$ & $1.95460789577555396208826644691$ \\
$125$ & $2.70242962410312278208767904044$ \\
$126$ & $2.09439510239319549230842892219$ \\
$127$ & $2.72368730766675765849410753531$ \\
$128$ & $1.78859193595207408964845885470$ \\
$129$ & $2.15300367350566872872181166589$ \\
$130$ & $2.12917795839603170191079901461$ \\
$131$ & $2.57848536120487307251550163471$ \\
$132$ & $1.76974436700717745842473969267$ \\
$133$ & $2.71715522600071392800943235026$ \\
$134$ & $2.40177874951129444901523775053$ \\
\hline
\end{tabular}
}
\scalebox{0.606}{
\begin{tabular}{|c|c|}
\hline
$q$  &  $M_q$\\ \hline 
$135$ & $2.11995352095278398486512754910$ \\
$136$ & $1.72431411431745665223798349732$ \\
$137$ & $2.72051117298997095021332225039$ \\
$138$ & $2.06213746115420018933849107128$ \\
$139$ & $2.78392266035840572927430273739$ \\
$140$ & $2.05172721872421812035646128332$ \\
$141$ & $2.29124192852861593669991478644$ \\
$142$ & $2.60986917715784586434512887899$ \\
$143$ & $2.62713175526143872638214799918$ \\
$144$ & $1.28254983016186409554403635967$ \\
$145$ & $2.43829876081692332338284132881$ \\
$146$ & $2.22099352575696724500691589494$ \\
$147$ & $2.02255384109558280368674411138$ \\
$148$ & $1.99550481523309639952094476513$ \\
$149$ & $2.68631183937556722239373025709$ \\
$150$ & $1.81177664489821466934352802615$ \\
$151$ & $2.62352855587439821150270052709$ \\
$152$ & $1.66331503401599345646761861198$ \\
$153$ & $2.36894896286376229372867859617$ \\
$154$ & $2.43523343388619715784941500753$ \\
$155$ & $2.68173130074102629208719570379$ \\
$156$ & $2.01222972650783293991026647857$ \\
$157$ & $2.91562362895517732081839286279$ \\
$158$ & $2.49865862902662621045751915333$ \\
$159$ & $2.49144503560985776378540431624$ \\
$160$ & $1.57079632679489661923132169164$ \\
$161$ & $2.70131898348000368760805411733$ \\
$162$ & $1.85460990010218743941830369297$ \\
$163$ & $2.69099736683125370993373921388$ \\
$164$ & $2.11431182971789691740711719452$ \\
$165$ & $2.19228536898779746002203139574$ \\
$166$ & $2.40279523907172221214735958325$ \\
$167$ & $2.74644085264695532114443109623$ \\
$168$ & $1.67925190836271397945223134714$ \\
$169$ & $2.52322934735426698644408474464$ \\
$170$ & $2.34387207448910709235492592796$ \\
$171$ & $2.41874068993618196904305168599$ \\
$172$ & $2.15300367350566872872181166589$ \\
$173$ & $2.83083393377236324187600452713$ \\
$174$ & $2.14059432707202598850826022652$ \\
$175$ & $2.71272588538338440992853033706$ \\
$176$ & $1.65533401256501131763033240531$ \\
$177$ & $2.37483432765382947109503698217$ \\
$178$ & $2.48752834330666367191078756009$ \\
$179$ & $2.95215347090837063371989208638$ \\
$180$ & $1.80945332629550203732956607281$ \\
$181$ & $2.55866549759635341623707612755$ \\
$182$ & $2.37482082344745189756925059590$ \\
$183$ & $2.28173329133748979800976062080$ \\
$184$ & $1.96520205410785916590276700512$ \\
$185$ & $2.79298789317197424525982798487$ \\
$186$ & $1.97538604461783624686153124312$ \\
$187$ & $2.79411040714987147526776186875$ \\
$188$ & $2.29124192852861593669991478644$ \\
$189$ & $2.32471019193358139074467540535$ \\
$190$ & $2.57856484292230594522273621834$ \\
$191$ & $2.95512966360404799352636309788$ \\
$192$ & $1.45122657606971550702884393708$ \\
$193$ & $2.60255291569166233786515416685$ \\
$194$ & $2.42085614235917869433795064647$ \\
$195$ & $2.29429488381819057673465872290$ \\
$196$ & $2.01688417594425696685721289539$ \\
$197$ & $2.81468933588096728324501080140$ \\
$198$ & $2.25006635803238006730577318676$ \\
$199$ & $2.79249308566493928174043020396$ \\
$200$ & $1.70951180694329794837835134968$ \\
\hline
\end{tabular}
}
\caption{\label{table1}
{\small
Values of $M_q$ for every  $3 \le q \le \boundtables$ with 
$30$-digit precision (the last printed digit is rounded by PARI/GP); computed with PARI/GP, v.~2.13.0.
Total computation time:   23 sec.,   471 millisecs.  
}
}
\end{table} 
\newpage

\begin{table}[H]
\scalebox{0.575}{
\begin{tabular}{|c|c|}
\hline
$q$  &  $m_q$\\ \hline 
$3$ & $0.604599788078072616864692752547$ \\
$4$ & $0.785398163397448309615660845820$ \\
$5$ & $0.430408940964004038889433232951$ \\
$6$ & $0.604599788078072616864692752547$ \\
$7$ & $0.547959686797993973084485988763$ \\
$8$ & $0.623225240140230513394020080251$ \\
$9$ & $0.604599788078072616864692752547$ \\
$10$ & $0.430408940964004038889433232951$ \\
$11$ & $0.618351934876807874060419662662$ \\
$12$ & $0.604599788078072616864692752547$ \\
$13$ & $0.598987497945465758207499250242$ \\
$14$ & $0.547959686797993973084485988763$ \\
$15$ & $0.430408940964004038889433232951$ \\
$16$ & $0.623225240140230513394020080251$ \\
$17$ & $0.453546340908733659287599108349$ \\
$18$ & $0.604599788078072616864692752547$ \\
$19$ & $0.413193436540565451244291268589$ \\
$20$ & $0.430408940964004038889433232951$ \\
$21$ & $0.547959686797993973084485988763$ \\
$22$ & $0.618351934876807874060419662662$ \\
$23$ & $0.552304916713058385866569568830$ \\
$24$ & $0.604599788078072616864692752547$ \\
$25$ & $0.430408940964004038889433232951$ \\
$26$ & $0.598987497945465758207499250242$ \\
$27$ & $0.604599788078072616864692752547$ \\
$28$ & $0.547959686797993973084485988763$ \\
$29$ & $0.451093787499735920935796126723$ \\
$30$ & $0.430408940964004038889433232951$ \\
$31$ & $0.440122223433962808617040155019$ \\
$32$ & $0.601117729884346270714091624669$ \\
$33$ & $0.604599788078072616864692752547$ \\
$34$ & $0.453546340908733659287599108349$ \\
$35$ & $0.430408940964004038889433232951$ \\
$36$ & $0.604599788078072616864692752547$ \\
$37$ & $0.420012687371836528710987165573$ \\
$38$ & $0.413193436540565451244291268589$ \\
$39$ & $0.598987497945465758207499250242$ \\
$40$ & $0.430408940964004038889433232951$ \\
$41$ & $0.531058786094205975539953409040$ \\
$42$ & $0.547959686797993973084485988763$ \\
$43$ & $0.479088388239857211764493892920$ \\
$44$ & $0.585417754933143702267998692855$ \\
$45$ & $0.430408940964004038889433232951$ \\
$46$ & $0.552304916713058385866569568830$ \\
$47$ & $0.367129807516487530311860512131$ \\
$48$ & $0.604599788078072616864692752547$ \\
$49$ & $0.498129254518742697997375454169$ \\
$50$ & $0.430408940964004038889433232951$ \\
$51$ & $0.453546340908733659287599108349$ \\
$52$ & $0.598987497945465758207499250242$ \\
$53$ & $0.413967522107356708589184273274$ \\
$54$ & $0.604599788078072616864692752547$ \\
$55$ & $0.430408940964004038889433232951$ \\
$56$ & $0.547959686797993973084485988763$ \\
$57$ & $0.413193436540565451244291268589$ \\
$58$ & $0.451093787499735920935796126723$ \\
$59$ & $0.332420580251333196195200170353$ \\
$60$ & $0.430408940964004038889433232951$ \\
$61$ & $0.365767556044524607545195023959$ \\
$62$ & $0.440122223433962808617040155019$ \\
$63$ & $0.547959686797993973084485988763$ \\
$64$ & $0.601117729884346270714091624669$ \\
$65$ & $0.430408940964004038889433232951$ \\
$66$ & $0.604599788078072616864692752547$ \\
$67$ & $0.383806628882915516388036615854$ \\
$68$ & $0.453546340908733659287599108349$ \\
\hline
\end{tabular}
}
\scalebox{0.575}{
\begin{tabular}{|c|c|}
\hline
$q$  &  $m_q$\\ \hline
$69$ & $0.552304916713058385866569568830$ \\
$70$ & $0.430408940964004038889433232951$ \\
$71$ & $0.445439715279504151396753278541$ \\
$72$ & $0.604599788078072616864692752547$ \\
$73$ & $0.332816572422086433683238786444$ \\
$74$ & $0.420012687371836528710987165573$ \\
$75$ & $0.430408940964004038889433232951$ \\
$76$ & $0.413193436540565451244291268589$ \\
$77$ & $0.421073984120137407835557128276$ \\
$78$ & $0.598987497945465758207499250242$ \\
$79$ & $0.428163851610403317805757592040$ \\
$80$ & $0.430408940964004038889433232951$ \\
$81$ & $0.482969705009469994409850936537$ \\
$82$ & $0.531058786094205975539953409040$ \\
$83$ & $0.462387538865549865563618907321$ \\
$84$ & $0.547959686797993973084485988763$ \\
$85$ & $0.430408940964004038889433232951$ \\
$86$ & $0.479088388239857211764493892920$ \\
$87$ & $0.451093787499735920935796126723$ \\
$88$ & $0.585417754933143702267998692855$ \\
$89$ & $0.389534102336872450885689912904$ \\
$90$ & $0.430408940964004038889433232951$ \\
$91$ & $0.446138400444371209202901940259$ \\
$92$ & $0.534749081131232483346513613662$ \\
$93$ & $0.440122223433962808617040155019$ \\
$94$ & $0.367129807516487530311860512131$ \\
$95$ & $0.413193436540565451244291268589$ \\
$96$ & $0.601117729884346270714091624669$ \\
$97$ & $0.389092591237496211100095150861$ \\
$98$ & $0.498129254518742697997375454169$ \\
$99$ & $0.556336740505493234121551435711$ \\
$100$ & $0.430408940964004038889433232951$ \\
$101$ & $0.356508197108894602900482691212$ \\
$102$ & $0.453546340908733659287599108349$ \\
$103$ & $0.397853727338476109339913228255$ \\
$104$ & $0.598987497945465758207499250242$ \\
$105$ & $0.430408940964004038889433232951$ \\
$106$ & $0.413967522107356708589184273274$ \\
$107$ & $0.404936054291382959432040324919$ \\
$108$ & $0.604599788078072616864692752547$ \\
$109$ & $0.395439291752596225073889921584$ \\
$110$ & $0.430408940964004038889433232951$ \\
$111$ & $0.420012687371836528710987165573$ \\
$112$ & $0.547959686797993973084485988763$ \\
$113$ & $0.429318876983891967381982481536$ \\
$114$ & $0.413193436540565451244291268589$ \\
$115$ & $0.430408940964004038889433232951$ \\
$116$ & $0.451093787499735920935796126723$ \\
$117$ & $0.520803376341663832109816800257$ \\
$118$ & $0.332420580251333196195200170353$ \\
$119$ & $0.412066798755779847198488618713$ \\
$120$ & $0.430408940964004038889433232951$ \\
$121$ & $0.434156735952171172172582423965$ \\
$122$ & $0.365767556044524607545195023959$ \\
$123$ & $0.449402261484866185812132591013$ \\
$124$ & $0.440122223433962808617040155019$ \\
$125$ & $0.430408940964004038889433232951$ \\
$126$ & $0.547959686797993973084485988763$ \\
$127$ & $0.370583958563266768017674320198$ \\
$128$ & $0.601117729884346270714091624669$ \\
$129$ & $0.479088388239857211764493892920$ \\
$130$ & $0.430408940964004038889433232951$ \\
$131$ & $0.352086298168602566529507073368$ \\
$132$ & $0.585417754933143702267998692855$ \\
$133$ & $0.395478409042985251655275470376$ \\
$134$ & $0.383806628882915516388036615854$ \\
\hline
\end{tabular}
}
\scalebox{0.575}{
\begin{tabular}{|c|c|}
\hline
$q$  &  $m_q$\\ \hline 
$135$ & $0.430408940964004038889433232951$ \\
$136$ & $0.453546340908733659287599108349$ \\
$137$ & $0.417505824794320762685433840755$ \\
$138$ & $0.552304916713058385866569568830$ \\
$139$ & $0.416862458745916605478311229554$ \\
$140$ & $0.430408940964004038889433232951$ \\
$141$ & $0.367129807516487530311860512131$ \\
$142$ & $0.445439715279504151396753278541$ \\
$143$ & $0.420201642019180957031600673770$ \\
$144$ & $0.604599788078072616864692752547$ \\
$145$ & $0.430408940964004038889433232951$ \\
$146$ & $0.332816572422086433683238786444$ \\
$147$ & $0.498129254518742697997375454169$ \\
$148$ & $0.420012687371836528710987165573$ \\
$149$ & $0.383901386511619291505439140039$ \\
$150$ & $0.430408940964004038889433232951$ \\
$151$ & $0.375982979158935926530947639675$ \\
$152$ & $0.413193436540565451244291268589$ \\
$153$ & $0.453546340908733659287599108349$ \\
$154$ & $0.421073984120137407835557128276$ \\
$155$ & $0.424744511168658171711726795906$ \\
$156$ & $0.598987497945465758207499250242$ \\
$157$ & $0.369949970311229129154525082455$ \\
$158$ & $0.428163851610403317805757592040$ \\
$159$ & $0.413967522107356708589184273274$ \\
$160$ & $0.430408940964004038889433232951$ \\
$161$ & $0.408568710916191166131009770920$ \\
$162$ & $0.482969705009469994409850936537$ \\
$163$ & $0.246068527552960243897853273760$ \\
$164$ & $0.498220344738919466049222615474$ \\
$165$ & $0.430408940964004038889433232951$ \\
$166$ & $0.462387538865549865563618907321$ \\
$167$ & $0.350631724723697517489360519141$ \\
$168$ & $0.547959686797993973084485988763$ \\
$169$ & $0.372013285127068443884524331707$ \\
$170$ & $0.430408940964004038889433232951$ \\
$171$ & $0.413193436540565451244291268589$ \\
$172$ & $0.479088388239857211764493892920$ \\
$173$ & $0.382663127794428634867406502671$ \\
$174$ & $0.451093787499735920935796126723$ \\
$175$ & $0.430408940964004038889433232951$ \\
$176$ & $0.503296496656052465507619897240$ \\
$177$ & $0.332420580251333196195200170353$ \\
$178$ & $0.389534102336872450885689912904$ \\
$179$ & $0.341096922920489890524209012663$ \\
$180$ & $0.430408940964004038889433232951$ \\
$181$ & $0.400130740444375343545002983965$ \\
$182$ & $0.446138400444371209202901940259$ \\
$183$ & $0.365767556044524607545195023959$ \\
$184$ & $0.534749081131232483346513613662$ \\
$185$ & $0.420012687371836528710987165573$ \\
$186$ & $0.440122223433962808617040155019$ \\
$187$ & $0.414644064612111527601828504332$ \\
$188$ & $0.367129807516487530311860512131$ \\
$189$ & $0.533350480933598026213622422760$ \\
$190$ & $0.413193436540565451244291268589$ \\
$191$ & $0.369946585080866896719707254253$ \\
$192$ & $0.601117729884346270714091624669$ \\
$193$ & $0.362378327133943195753957944776$ \\
$194$ & $0.389092591237496211100095150861$ \\
$195$ & $0.430408940964004038889433232951$ \\
$196$ & $0.498129254518742697997375454169$ \\
$197$ & $0.373837653864628774889214137854$ \\
$198$ & $0.556336740505493234121551435711$ \\
$199$ & $0.407431952492485762271173360682$ \\
$200$ & $0.430408940964004038889433232951$ \\

\hline
\end{tabular}
} 
\caption{\label{table2}
{\small
Values of $m_q$ for every  $3 \le q \le \boundtables$ with 
$30$-digit precision (the last printed digit is rounded by PARI/GP); computed with PARI/GP, v.~2.13.0.
Total computation time:   23  sec., 696 millisecs.  
}
}
\end{table}  

  \vfill\eject

  \end{document}